\documentclass[12pt]{amsart}

\usepackage{amsmath}
\usepackage{amssymb}
\usepackage{amsthm}
\usepackage{amscd}
\usepackage[all]{xy}
\usepackage[utf8]{inputenc}
\usepackage[T1]{fontenc}   

\usepackage{color}
\usepackage{hyperref}

\usepackage{mathpazo}
\usepackage[scaled=.95]{helvet}
\usepackage{courier}

\swapnumbers
\headsep=1cm \numberwithin{equation}{section}
\newtheorem{theorem}{Theorem}[section]
\newtheorem{conjecture}[theorem]{Conjecture}
\newtheorem{proposition}[theorem]{Proposition}
\newtheorem{lemma}[theorem]{Lemma}
\newtheorem{question}[theorem]{Question}
\newtheorem{corollary}[theorem]{Corollary}

\theoremstyle{definition}
\newtheorem{definition}[theorem]{Definition}
\newtheorem{example}[theorem]{Example}

\theoremstyle{remark}
\newtheorem{remark}[theorem]{Remark}

\usepackage[all]{xy}

\usepackage{amsmath}
\usepackage{amssymb}
\usepackage{amsthm}
\usepackage{amscd}
\usepackage{color}
\usepackage{hyperref}

\newcommand\Hom{\operatorname{Hom}}
\newcommand\Ext{\operatorname{Ext}}
\newcommand\Tor{\operatorname{Tor}}

\newcommand\depth{\operatorname{depth}}

\newcommand\Ker{\operatorname{\Ker}}

\numberwithin{equation}{section}

     \title[The ARC, Finite $C$-injective dimension, and Vanishing of $\Ext$]{The Auslander-Reiten Conjecture, Finite $C$-Injective Dimension of $\operatorname{Hom}$, and Vanishing of $\Ext$}

     \author{Victor D. Mendoza-Rubio}
     \address{Universidade de S{\~a}o Paulo -
ICMC, Caixa Postal 668, 13560-970, S{\~a}o Carlos-SP, Brazil}
     \email{vicdamenru@usp.br}
     \thanks{The first author was supported by grant 2022/03372-5, São Paulo Research Foundation (FAPESP)}

     \author{Victor H. Jorge-Pérez}
     \address{Universidade de S{\~a}o Paulo -
ICMC, Caixa Postal 668, 13560-970, S{\~a}o Carlos-SP, Brazil}
     \email{vhjperez@icmc.usp.br}
     \thanks{The second  author was supported by grant 2019/21181-0, São Paulo Research Foundation (FAPESP) }

     \date{Month, Day, Year}

     \keywords{Auslander-Reiten conjecture, vanishing of Ext, injective dimension, semidualizing module, $C$-injective dimension,  canonical module, Gorenstein ring}
     \subjclass[2020]{Primary: 13D45, 13D07, 13C10, 13C14; Secondary: 13D05, 13D02, 13H10, 14B15.}

\begin{document}
     \begin{abstract}
     Let $R$ be a Noetherian local ring,  and let $C$ be a semidualizing $R$-module. In this paper, we present some results concerning the vanishing of $\operatorname{Ext}$ and finite injective dimension of $\operatorname{Hom}$. Additionally, we extend these results in terms of finite  $C$-injective dimension of $\operatorname{Hom}$. We also investigate the consequences of some of these extensions in the case where $R$ is Cohen-Macaulay and $C$ is a canonical module for $R.$ Furthermore, we provide positive answers to the Auslander-Reiten conjecture for finitely generated $R$-modules $M$ such that $\mathcal{I}_C\operatorname{-id}_R(\Hom_R(M,R))<\infty$ or $M \in \mathcal{A}_C(R)$ with $\mathcal{I}_C \operatorname{-id}_R(\operatorname{Hom}_R(M,M))<\infty$. Moreover, we derive a number of criteria for a semidualizing $R$-module $C$ to be a canonical module for $R$ in terms of the vanishing of $\operatorname{Ext}$ and the finite $C$-injective dimension of $\Hom$.
     \end{abstract}
     \maketitle

\section{Introduction}
Throughout this paper, we assume $R$ to be a commutative Noetherian local ring of dimension $d$ with maximal ideal $\mathfrak{m}$. Moreover, all $R$-modules are considered to be finitely generated, unless otherwise stated. 

The study of the vanishing of $\Ext$ modules over a local ring is an actively researched topic in commutative algebra and homological algebra. One of the most studied and important conjectures in commutative algebra is the  Auslander-Reiten conjecture:

\begin{conjecture}\cite{OnageneralizedversionoftheNakayamaconjecture}\label{AR}
If $M$ is an $R$-module such that $\Ext_R^i(M,R)=\Ext_R^i(M,M)=0$ for every $i\geq 1$, then $M$ is free.
\end{conjecture}

So far, many criteria for a given module to be free have been described in terms of the vanishing of $\Ext$ modules. In fact, this conjecture is still open but has been extensively studied and resolved in several cases.   References of some of them can be found in \cite{AuslanderReitenConjectureforNormalRings, injectivedimensiontakahashi}.

Recently, Ghosh and Takahashi, in \cite{injectivedimensiontakahashi}, provided criteria for an $R$-module to be free in terms of vanishing of $\operatorname{Ext}$ and finite injective dimension of $\operatorname{Hom}$, and proved that the Auslander-Reiten conjecture holds true when at least one of the modules $\Hom_R(M,R)$ and $\Hom_R(M,M)$ has finite injective dimension. Moreover, motivated by \cite[Theorem 2.5]{injectivedimensiontakahashi}, they posed the following intriguing question:
\begin{question}\cite[Question 2.9]{injectivedimensiontakahashi}\label{questio2.9InjeFTakahashi}
Let $M$ and $N$ be nonzero modules over $R$. If $\Hom_R(M,N)$ has finite injective dimension and $\Ext_R^i(M,N)=0$ for every $1\leq i \leq d$, then is $M$ free and does $N$ have finite injective dimension?
\end{question}
Inspired by the Auslander-Reiten Conjecture and Question \ref{questio2.9InjeFTakahashi}, as well as the results obtained by Ghosh and Takahashi in \cite[Theorems 2.5 and 2.15]{injectivedimensiontakahashi}, in this work, we aim to improve and extend these results in terms of the finite $C$-injective dimension of $\Hom$, where $C$ is a semidualizing $R$-module. Consequently, we provide new cases where the Auslander-Reiten conjecture is true. More explicitly, we prove the following results. 
\begin{theorem}[=Theorem \ref{otrageneralizacionde2.15}] \label{teo2.2}Let $M$ and $N$ be nonzero $R$-modules and let $t=\depth(N)$. Suppose that $\Hom_R(M,N)$  has finite injective dimension, and that  $\Ext_R^i(M,N)=\Ext_R^j(M,R)=0$ for all $1\leq i \leq t$ and  $1\leq j \leq d$. Then $M$ is free and $N$ has finite injective dimension. 
\end{theorem}
 \begin{theorem}[=Theorem \ref{estasigeneralizacionde215}]  \label{mejorade2.15}
Let  $M$ be a nonzero $R$-module such that  $\Hom_R(M,M)$ has finite injective dimension, and $\Ext_R^i(M,M)=\Ext_R^j(M,R)=0$
for all $1\leq i \leq d-1$ and $1\leq j \leq d$. Then $M$ is free and $R$ is Gorenstein.
\end{theorem}
\begin{theorem}[=Theorem \ref{generalizaciondoteo2.5}] \label{extensionde2.5}
     Let $M$ and $N$ be nonzero $R$-modules such that $N \in \mathcal{A}_C(R)$ is maximal Cohen-Macaulay. Suppose that $\Hom_R(M,N)$ has finite $C$-injective dimension and that $\Ext_R^i(M,N)=0$ for all $1\leq i \leq d$. Then $R$ is Cohen-Macaulay, $M$ is free and $N$ has finite $C$-injective dimension. 
 \end{theorem}
 \begin{theorem}[=Theorem \ref{2.15generalizadoIC}] \label{extensionymejorade2.15}
Let $M$ be a nonzero $R$-module such that $M \in \mathcal{A}_C(R)$.  Suppose that
 $\Hom_R(M,M)$ has finite $C$-injective dimension, and that $\Ext_R^i(M,M)=\Ext_R^j(M,R)=0$ for all $1\leq i \leq d-1$ and $1\leq j \leq d$.  Then  $R$ is Cohen-Macaulay,  $M$ is free and $C$ is a canonical module for $R$. 
\end{theorem}
As consequences of the theorems mentioned above, we not only improve, extend, or generalize results presented in \cite{injectivedimensiontakahashi}, but we also offer partial answers to Question \ref{questio2.9InjeFTakahashi} and the Auslander-Reiten conjecture. In particular, let us highlight the following outcomes. Theorem \ref{teo2.2} provides a generalization of \cite[Corollary 2.10(2)]{injectivedimensiontakahashi}. Theorem \ref{mejorade2.15} represents an improvement over \cite[Theorem 2.15]{injectivedimensiontakahashi}. Theorems \ref{extensionde2.5} and \ref{extensionymejorade2.15} extend or enhance the results in \cite[Theorems 2.5 and 2.15]{injectivedimensiontakahashi}, taking into account the $C$-injective dimension of $\Hom$. Partial answers to Question \ref{questio2.9InjeFTakahashi} are presented in Corollaries \ref{Exemplos} and \ref{corconNrigid}, while partial answers to the Auslander-Reiten conjecture are found in Corollaries \ref{ARCparaCidimHom(M,R)}, \ref{ARCpdHom(M,R)}, \ref{ARCforCidHom(M,M)}, and \ref{ARGorenstein}. Furthermore, we establish criteria for an $R$-module to be free, based on the vanishing of $\Ext$ and the finite projective dimension of $\Hom$. This is elaborated in Corollaries \ref{corpdGprojective} and \ref{otrocritpd}, as well as Theorem \ref{pdHom}. Notably, Theorem \ref{pdHom} generalizes \cite[Corollary 2.14]{injectivedimensiontakahashi}. Additionally, we provide criteria for a semidualizing module to be considered a canonical module of a Cohen-Macaulay local ring,  offering an improved version of \cite[Theorem 3.6]{injectivedimensiontakahashi}. These criteria and this improved version are established in Theorem \ref{caracterizacionessemidualizingserdualizing} and Corollary \ref{Gorentein}, respectively.

The organization of this paper is as follows: In Section 2, we introduce the notation, definitions, and some necessary known results that will be used in this paper.

In Section 3, we investigate the consequences of imposing finite injective dimension on $\operatorname{Hom}_R(M,N)$ under the conditions $\operatorname{Ext}_R^i(M,N)=0$ for $i=1,\dots,n$ for some $n\geq 1$ and $\operatorname{Ext}_R^i(M,R)=0$ for all $i=1,\dots,d$. As an application of the results obtained, we improve the result established by Ghosh and Takahashi in \cite[Theorem 2.15]{injectivedimensiontakahashi}. Additionally, we provide some partial positive answers to Question \ref{questio2.9InjeFTakahashi}, for instance in the case where $N$ is Tor-rigid.

In Section 4, we provide an overview of semidualizing modules and $C$-injective dimension, where $C$ is a semidualizing module. 

In Section 5, we employ the concept of semidualizing modules to extend certain results obtained in Section 
3 and \cite{injectivedimensiontakahashi} to the context of finite $C$-injective dimension of $\operatorname{Hom}$, for a semidualizing $R$-module $C$. Additionally, as applications of the extensions obtained,  we study the consequences of some of them in the case where $R$ is Cohen-Macaulay and $C$ is a canonical module for $R$, and provide partial positive answers to the Auslander-Reiten conjecture. In the last section, we provide criteria for a semidualizing $R$-module $C$ to be considered a canonical module for $R$ in terms of finite $C$-injective dimension of certain Hom.  Additionally, we present an improved version of \cite[Theorem 3.6]{injectivedimensiontakahashi}.

\begin{remark}
    The first version of this paper was finished and submitted to a journal in June 2023. Later, in October  2023, Ghosh and Dey submitted to arXiv the work \cite{FiniteHomologicalDimensionOfHomv2}, where we noticed that they independently obtained some similar results to ours. 
\end{remark}

\section{Setup and Background}

In this section, we provide essential definitions and properties that are used in this paper.

Here, the notation $M^\ast$ represents the algebraic dual of $M$, which is defined as $\Hom_R(M,R)$.  Let  $M$ be an $R$-module and consider a minimal free resolution $$\cdots \longrightarrow F_i \stackrel{\varphi_i}{\longrightarrow} F_{i-1} \longrightarrow \cdots \longrightarrow F_1 \stackrel{\varphi_1}{\longrightarrow} F_0 \stackrel{\varphi_0}{\longrightarrow} M \longrightarrow 0$$
of $M$.  For $i\geq 1$, the $i$-\textit{syzygy} of $M$, denoted by $\Omega^i(M)$, is defined as the kernel of the map $\varphi_{i-1}$. When $i=0$, we set $\Omega^0(M)=M$. For $i\geq 0$, the modules $\Omega^i(M)$ are defined uniquely up to isomorphism. If $N$ is an $R$-module, we say that $M$ and $N$ are \textit{stably isomorphic} and write $M\approx N$  if there exist free $R$-modules $F$ and $G$ such that $M \oplus F \cong N \oplus G$. For $i\geq 0$, we say that $M$ is an $i$-syzygy if $M \approx \Omega^i(N)$ for some $R$-module $N$. 

The \textit{Auslander Transpose} of $M$, denoted by $\operatorname{Tr}(M)$, is defined as the cokernel of  the induced map $\varphi_{1}^\ast: F_0^\ast \to F_1^\ast.$ The $R$-module $\operatorname{Tr}(M)$ is defined uniquely up to isomorphism. It is easy to see that $\operatorname{Tr}(\operatorname{Tr}(M))\approx M$. For $i\geq 1$, set  $\mathcal{T}_i(M)=\operatorname{Tr}(\Omega^{i-1}(M))$. Then for   each  $i\geq 0$ and $R$-module $N$, there exists an exact sequence \cite[Theorem 2.8]{Stablemoduletheory}
\begin{equation}\label{auslanderBridgerSequence}
\operatorname{Tor}_2^R\left(\mathcal{T}_{i+1} (M),N \right) \rightarrow \operatorname{Ext}_R^i(M, R) \otimes_RN \rightarrow \operatorname{Ext}_R^i(M,N) \rightarrow \operatorname{Tor}_1^R\left(\mathcal{T}_{i+1} (M),N\right) \rightarrow 0.
\end{equation}

The notion of Gorenstein dimension was introduced by Auslander \cite{AnneauxdeGorensteinettorsionenalgebrecommutative}  and developed by Auslander and Bridger in \cite{Stablemoduletheory}. An $R$-module $M$ is said to be $G$-\textit{projective} if the natural map $M \to M^{ \ast \ast}$ is an isomorphism and $\Ext_R^i(M,R)=\Ext_R^i(M^\ast,R)=0$ for all $i\geq 1$. The \textit{Gorenstein dimension} of $M$, denoted by $\operatorname{G-dim}_R(M)$, is defined to be the infimum of all nonnegative integers $n$, such that there exists an exact sequence
$$
0 \rightarrow G_n \rightarrow \cdots \rightarrow G_0 \rightarrow M \rightarrow 0,
$$
where each $G_i$ is $G$-projective. It is easy to see that  $\operatorname{G-dim}_R(M)=0$ if and only if $M$ is $G$-projective. If $M\not=0$ and $\operatorname{G-dim}_R(M)<\infty$, then $\operatorname{G-dim}_R(M)=\operatorname{depth}(R)-\operatorname{depth}(M)$ (see \cite[Theorem 1.4.8]{GorensteinDimensions}) and this equality is known as the \textit{Auslander-Bridger formula}. If $R$ is Gorenstein, then every $R$-module has finite $G$-dimension (see \cite[Theorem 1.4.9]{GorensteinDimensions}). 

An $R$-module $N$ is \textit{Tor-rigid} if for every $R$-module $M$ and every $i\geq 1$ holds the following implication: $$\Tor_{i}^R(M,N)=0\implies \Tor_{i+1}^R(M,N)=0.$$
For examples of Tor-rigid modules, we refer the reader to \cite[Example 3.4]{VanishingOfCohomologyOverCompleteIntersectionRings} and \cite{SomeHomolocalPropertiesofModulesoveraCompleteIntersectionwithApplications}.

The following result will be used frequently in this paper.
\begin{theorem}[Bass' Theorem]
\label{Bassconjecture}
If $M$ is a nonzero $R$-module of finite injective dimension, then  $R$ is Cohen-Macaulay and $\operatorname{id}_R(M)=d$.  
\end{theorem}
\begin{proof}
    This follows from \cite[Theorem 1.6]{injectivedimensiontakahashi} and  \cite[Theorem 3.1.17]{bruns}.
\end{proof}

\section{Finite injective dimension of $\Hom$ and vanishing of $\Ext$ }
In this section, for two $R$-modules $M$ and $N$, we explore the consequences of $\Hom_R(M,N)$ having finite injective dimension under the condition that $\Ext_R^i(M,N)=0$ for $i=1,\dots,n$ for some $n\geq 1$, and $\Ext_R^i(M,R)=0$ for $i=1,\dots,d$. Additionally, we explore Question \ref{questio2.9InjeFTakahashi}. For this, we need the following lemma, where items (1) and (2) are probably well-known to the experts and used repeatedly in the literature, while item (3) is a recent result of \cite{MCMtensorproductsandvanishingofExtmodules}.

\begin{lemma}\label{dephdeHom}
    Let $M$ and $N$ be nonzero $R$-modules.
    \begin{enumerate}
     \item[(1)] Let $\pmb{x}=x_1,\ldots,x_s$ be an $N$-sequence.  If  $\Hom_R(M,N)\not=0$ and  $\Ext_R^i(M,N)=0$ for all $1\leq i \leq s$, then $\pmb{x}$ is a $\Hom_R(M,N)$-sequence and 
    $$\Hom_R(M,N/\pmb{x}N) \cong \frac{\Hom_R(M, N)}{\pmb{x}\Hom_R(M,N)}.$$
      \item[(2)] If $\Ext_R^i(M,N)=0$ for all $1\leq i \leq \depth(N)$, then $\Hom_R(M,N)\not=0$ and $$\depth(N)=\depth(\Hom_R(M,N)).$$
    \item[(3)]  If $\depth(N)=0$ and $\Ext_R^1( \operatorname{Tr}(M), \Hom_R(M,N))=0$, then $M$ is free. 
    \end{enumerate}
\end{lemma}
\begin{proof}
$(1)$  By induction,  it is enough to consider $s=1$. Since $x_1$ is an $N$-regular element, then the sequence  $$\xymatrix{0 \ar[r]& N \ar[r]^{x_1}& N \ar[r]& N/x_1N \ar[r]&0}  $$
is exact. Once $\Ext_R^1(M,N)=0$, this sequence induces an exact sequence 
\begin{equation*}
    \xymatrix{0 \ar[r]& \Hom_R(M,N) \ar[r]^{x_1}& \Hom_R(M,N) \ar[r]& \Hom_R(M,N/x_1N)\ar[r]& 0   }.
    \end{equation*}
    Thus, $x_1$ is a $\Hom_R(M,N)$-regular element and 
$$\frac{\Hom_R(M,N)}{x_1\Hom_R(M,N)} \cong \Hom_R(M,N/x_1N).$$

$(2)$ Let $t=\depth(N)$.  Set $I=\operatorname{ann}(M)$.  Note that $I$ is a proper ideal of $R$ since $M\neq 0$. Since $N\neq 0$, by Nakayama's lemma, we have $IN\neq N$. Thus, $\operatorname{grade}(I, N)=\inf\{i\geq 0: \Ext_R^i(M,N)\neq 0\}$, and $0 \leq \operatorname{grade}(I,N) \leq t$. As $\Ext_R^i(M,N)=0$ for all $1\leq i \leq t$, it follows that $\operatorname{grade}(I,N)=0$, which means $\Hom_R(M,N)\neq 0.$ 

Now, as $\depth(N)=t$, then there exists an $N$-sequence $\pmb{x}=x_1, \ldots,x_t$. By (1),   $\pmb{x}$ is a $\Hom_R(M,N)$-sequence and $$\Hom_R(M,N/\pmb{x}N) \cong \frac{\Hom_R(M, N)}{\pmb{x}\Hom_R(M,N)}.$$
Hence, $$\depth(\Hom_R(M, N/\pmb{x}N))=\depth(\Hom_R(M,N))-t.$$
Since $\depth(N/\pmb{x}N)=0$, then $\depth( \Hom_R(M, N/\pmb{x}N))=0$ (\cite[Exercise 1.2.27]{bruns}). Therefore, \linebreak $\depth(\Hom_R(M,N))=t$.

 $(3)$ It follows from \cite[Proposition 3.3(2)]{MCMtensorproductsandvanishingofExtmodules}. 

\end{proof}

Note that Lemma \ref{dephdeHom} improves the given remark in \cite[Remark 2.2]{injectivedimensiontakahashi}. With this, we derive the following theorem, which generalizes a result of Ghosh and Takahashi (\cite[Corollary 2.10(2)]{injectivedimensiontakahashi}).
\begin{theorem}\label{otrageneralizacionde2.15} Let $M$ and $N$ be nonzero $R$-modules, and let $t=\depth(N)$. Suppose that $\Hom_R(M,N)$  has finite injective dimension, and that  $\Ext_R^i(M,N)=\Ext_R^j(M,R)=0$ for all $1\leq i \leq t$ and  $1\leq j \leq d$.  Then $M$ is free and $N$ has finite injective dimension. 
\end{theorem}
\begin{proof} 
By Lemma \ref{dephdeHom}(2), we have $\Hom_R(M,N)\neq 0$. By  Bass' theorem, we know that $R$ is Cohen-Macaulay and $\operatorname{id}_R (\Hom_R(M,N))=d$. Now, let $\pmb{x}=x_1,\dots,x_t$ be an $R$- and $N$-sequence. By Lemma \ref{dephdeHom}(1), we get 
    $$\Hom_R(M,N/\pmb{x}N)\cong \Hom_R(M,N)/\pmb{x}\Hom_R(M,N)$$
and $\pmb{x}$ is a $\Hom_R(M,N)$-sequence. This implies that $\Hom_R(M, N/\pmb{x}N)$ has finite injective dimension according to \cite[Exercise 4.3.3]{weibel}.  Now, consider a minimal free resolution of $M$:
$$\cdots \longrightarrow F_i \stackrel{\varphi_i}{\longrightarrow} F_{i-1} \longrightarrow \cdots \longrightarrow F_1 \stackrel{\varphi_1}{\longrightarrow} F_0 \stackrel{\varphi_0}{\longrightarrow} M \longrightarrow 0.$$ 
Since $\Ext_R^j(M,R)=0$ for all $1\leq j \leq d$, we have an exact sequence
    $$0  \longrightarrow M^\ast \stackrel{\varphi_0^\ast}{\longrightarrow} F_{0}^\ast \stackrel{\varphi_1^\ast}{\longrightarrow}  F_1^\ast \longrightarrow \cdots \stackrel{\varphi_{d+1}^\ast} \longrightarrow F_{d+1}^\ast \longrightarrow L  \longrightarrow 0,$$
where $L=\operatorname{coker}(\varphi_{d+1}^\ast)$. Since $\operatorname{Tr}(M)=\operatorname{coker}(\varphi_1^\ast)$, this exact sequence implies that $\operatorname{Tr}(M)$ is a $d$-syzygy of $L$, i.e., $\operatorname{Tr}(M) \approx \Omega^d L$. As $\operatorname{id}_R(\Hom_R(M,N/\pmb{x}N))=d$, we have $\Ext_R^{d+1}(L, \Hom_R(M,N/\pmb{x}N))=0$. But since $\operatorname{Tr}(M)\approx \Omega^d L$, we also have $\Ext_R^1\left( \operatorname{Tr}(M), \Hom_R(M,N/\pmb{x}N)\right)=0$. By Lemma \ref{dephdeHom}(3), this implies that $M$ is free. Consequently, $N$ has finite injective dimension since $\Hom_R(M,N)$ does.
\end{proof}

Motivated by the assumptions  of Theorem \ref{otrageneralizacionde2.15}, the following example presents  two cases in which an $R$-module $M$ satisfies $\Ext_R^i(M,R)=0$ for all $1\leq i\leq d$.

    \begin{example} \label{exemplosExt(M,R)=0} Let $M$ be an $R$-module. Then $\Ext_R^j(M,R)=0$ for $1\leq j \leq d$ in each one of the following situations:
    \begin{enumerate}
    \item[(1)] $R$ is Gorenstein and $M$ is maximal Cohen-Macaulay. 
    \item[(2)] $N$ is a nonzero Tor-rigid $R$-module such that $\Ext_R^i(M,N)=0$ for all $1\leq i \leq d$.
    \end{enumerate}
    \end{example}
\begin{proof}

$(1)$ By \cite[Exercise 3.1.24]{bruns}, $\Ext_R^i(M,R)=0$ for all $1\leq i \leq d$.

$(2)$ Suppose that $d\geq 1$. Fix an integer $j$ with $1\leq j \leq d$. Then $\Ext_R^j(M,N)=0$. By the exact sequence \eqref{auslanderBridgerSequence}, we have $\Tor_1^R\left( \mathcal{T}_{j+1}(M),N\right)=0$. Therefore, by the rigidity of $N$, we can conclude that $\Tor_2^R\left(\mathcal{T}_{j+1}(M),N\right)=0$. Consequently, from the exact sequence \eqref{auslanderBridgerSequence}, we have $\Ext_R^j(M,R) \otimes_R N=0$. Since $N\neq 0$, it follows that $\Ext_R^j(M,R)=0$. 
\end{proof}

\begin{corollary}\label{Exemplos}
    Let $R$ be a Gorenstein local ring, and let  $M$ and $N$ be nonzero $R$-modules such that $M$ is maximal Cohen-Macaulay. Let $t=\depth(N)$. Suppose that $\Hom_R(M,N)$ has finite injective dimension and that $\Ext_R^i(M,N)=0$ for all $1\leq i \leq t$.  Then $M$ is free and $N$ has finite injective dimension.
\end{corollary}
\begin{proof}
 This  follows from Theorem \ref{otrageneralizacionde2.15} and Example \ref{exemplosExt(M,R)=0}(1). 
\end{proof}
\begin{corollary}\label{corconNrigid}
    Let $M$ and $N$ be nonzero $R$-modules such that $N$ is Tor-rigid.  Suppose that $\Hom_R(M,N)$  has finite injective dimension and that  $\Ext_R^i(M,N)=0$ for all $1\leq i \leq d$. Then $M$ is free and $N$ has finite injective dimension. 
\end{corollary}
\begin{proof}
    This follows from Theorem \ref{otrageneralizacionde2.15} and Example \ref{exemplosExt(M,R)=0}(2). 
\end{proof}

It should be noted that Theorem \ref{otrageneralizacionde2.15} and its corollaries provide a partial answer to Question \ref{questio2.9InjeFTakahashi}.

Now, as an application of Theorem \ref{otrageneralizacionde2.15}, we can improve one of the main results given by Ghosh and Takahashi in \cite[Theorem 2.15]{injectivedimensiontakahashi} as follows.

\begin{theorem}\label{estasigeneralizacionde215}
Let  $M$ be a nonzero $R$-module such that $\Hom_R(M,M)$ has finite injective dimension, and $\Ext_R^i(M,M)=\Ext_R^j(M,R)=0$
for all $1\leq i \leq d-1$ and $1\leq j \leq d$. Then $M$ is free and $R$ is Gorenstein.
\end{theorem}
\begin{proof}
We may assume that $R$ is complete. By Bass' theorem, $R$ is Cohen-Macaulay, and hence admits a canonical module $\omega_R$. Let $t=\depth(M)$. We claim that $M$ is maximal Cohen-Macaulay, that is, $t=d$. Indeed, if $t\leq d-1$,  by hypothesis, $\Ext_R^i(M,M)=0$ for all $1\leq i \leq t$. Therefore, from Theorem \ref{otrageneralizacionde2.15}, $M$ is free, and hence $t=d$, which leads to a contradiction.

Now, by \cite[Theorem 2.3]{injectivedimensiontakahashi}, $M \cong \Gamma_\mathfrak{m}(M) \oplus R^r$ for some $r\geq 0$, and $M$ has finite injective dimension. But since $M$ is maximal Cohen-Macaulay, then $\Ext_R^d(M, \omega_R)=0$. Thus, the local duality theorem shows that $\Gamma_\mathfrak{m}(M)=0$. Hence, $M \cong R^r$. Since $M$ has finite injective dimension, it follows that $R$ is Gorenstein.
\end{proof}

 We can also observe that Theorem \ref{estasigeneralizacionde215} allows us to improve \cite[Corollary 2.14 and Theorem 3.6]{injectivedimensiontakahashi}.

 Recently, Zargar and Gheibi proved in \cite{NumericalApectsofComplexesOfFiniteHomologicalDimensions} the following result: If $M$ and $N$ are two $R$-modules such that $\operatorname{Ext}_R^i(M, N)=0$ for sufficiently large $i$ and $\operatorname{Ext}_R^i(M, N)$ has finite injective dimension for all $i$, then $\operatorname{pd}_R (M)$ and $\mathrm{id}_R(N)$ are finite. Using this result, we have the following proposition, which provides another partial answer to Question \ref{questio2.9InjeFTakahashi}. Consequently, Question \ref{questio2.9InjeFTakahashi} holds when $R$ is a local complete intersection of dimension $d\geq 1$ and codimension $d-1$.

\begin{proposition}
    Let $R$ be a local complete intersection of codimension $c$. Let $M$ and $N$ be nonzero $R$-modules. Suppose that:
    \begin{enumerate}
        \item  $\Hom_R(M,N)$ has finite injective dimension.
        \item $\Ext_R^i(M,N)=0$ for all $1\leq i \leq 1+c$.
       \item $\depth(N)\leq 1+c$.
    \end{enumerate}
    Then $M$ is free and $N$ has finite injective dimension.  
\end{proposition}
\begin{proof} 
From \cite[Theorem 4.1]{VanishingOfCohomologyOverCompleteIntersectionRings}, it follows that $\operatorname{Ext}_R^i(M,N)=0$ for all $i\geq 1$. Additionally, \cite[Corollary 4.3]{NumericalApectsofComplexesOfFiniteHomologicalDimensions} implies that $\operatorname{pd}_R(M)<\infty$ and $\operatorname{id}_R(N)<\infty$. Now, according to \cite[p. 154, Lemma 1(iii)]{CommutativeRingTheory}, we have $\operatorname{pd}_R(M)=\sup\{i\geq 0:\operatorname{Ext}_R^i(M,N)\neq 0\}$. Since $\operatorname{Ext}_R^i(M,N)=0$ for all $i\geq 1$, it follows that $\operatorname{pd}_R(M)=0$, which implies that $M$ is a free module.
\end{proof}

\section{Semidualizing modules and $C$-injective dimension}
In this section, we review some definitions and results concerning semidualizing modules and $C$-injective dimension, where $C$ is a semidualizing  $R$-module. This review enables us to generalize the results obtained in the previous section, as well as the main findings presented by Ghosh and Takahashi in \cite{injectivedimensiontakahashi}. It is worth noting that semidualizing modules were initially studied by Foxby \cite{GorensteinModulesAndRelatedModules} and Vasconcelos \cite{DivisorTheoryinModuleCategories}, and independently by Golod \cite{GdimensionandGeneralizedPerfectIdeals}. These modules play a crucial role in extending the concept of dualizing modules.

\begin{definition} An $R$-module $C$ is \textit{semidualizing} if the natural map $R\to \Hom_R(C,C)$ is an isomorphism and $\Ext_R^i(C,C)=0$ for all $i\geq 1$. 
\end{definition}
The ring $R$ is considered a semidualizing $R$-module. Moreover, if $R$ is a Cohen-Macaulay local ring with a canonical module $\omega_R$, then $\omega_R$ is a semidualizing $R$-module. It is worth noting that a canonical module of a Cohen-Macaulay local ring can be characterized as a semidualizing $R$-module of finite injective dimension.

For convenience, from now on, $C$ denotes a semidualizing $R$-module.

The following proposition contains some known results concerning semidualizing modules.
\begin{proposition}\cite{SemidualizingModules} \label{Propbasicasemidualizantes} The following statements hold:
\begin{enumerate}
    \item[(1)] $\depth(C)=\depth(R)$ and $\dim(C)=\dim(R)$. In particular, if $R$ is Cohen-Macaulay, then $C$ is maximal Cohen-Macaulay.
    \item[(2)] $C$ is indecomposable.
   \item[(3)] Let $M$ be an $R$-module. Then $M \neq 0$ if and only if $\operatorname{Hom}_R(C, M) \neq 0$. 
    \item[(3)] If $\varphi:R \to S$ is a flat ring homomorphism of Noetherian rings, then $C \otimes_R S$ is a semidualizing $S$-module. 
\end{enumerate}
\end{proposition}

\begin{definition}
    Let $M$ be an $R$-module. 
    \begin{enumerate}
    \item An $R$-module $M$ is called $C$-\textit{injective} if $M\cong \Hom_R(C,I)$ for some injective $R$-module $I$.
        \item A complex of the form $X_\bullet:0 \to M \to B^0 \to B^1 \to \cdots \to B^n\to \cdots,$
        where each $B^i$ is $C$-injective $R$-module, is called an $\mathcal{I}_C$-\textit{injective resolution}  if the complex $$C \otimes_R X_\bullet: 0 \to C \otimes_R M\to C \otimes_R B^0 \to C \otimes_R B^1 \to \cdots \to C \otimes_R B^n \to \cdots $$ is exact. The $\mathcal{I}_C$-injective dimension of $M$, denoted by $\mathcal{I}_C\operatorname{-id}_R(M)$, is defined as the infimum of all $n\geq 0$ for which there exists an $\mathcal{I}_C$-injective resolution  of the form $X_\bullet:0 \to M \to B^0 \to B^1 \to \cdots \to B^n\to 0.$
    \end{enumerate}
\end{definition}

Note that from Theorem \ref{formulaIcPc}, when $C=R$, we observe that the $C$-injective dimension coincides with injective dimension.

\begin{definition} The \textit{Auslander class} with respect to $C$, denoted by $\mathcal{A}_C(R)$, is the class of all (not necessarily finitely generated) $R$-modules $M$ such that:
\begin{enumerate}
    \item[(1)] The natural map $M \to \Hom_R(C, C \otimes_R M)$ is an isomorphism.
    \item[(2)] $\Tor_i^R(C,M)=0$ for all $i\geq 1.$
    \item[(3)] $\Ext_R^i(C, C \otimes_R M)=0$ for all $i\geq 1$.
\end{enumerate}
\end{definition}

If $C$ is a semidualizing $R$-module, note that $R \in \mathcal{A}_C(R)$.  When $R$ is Cohen-Macaulay and $C$ is a canonical module for $R$, then $\mathcal{A}_C(R)$ contains all $R$-modules of finite Gorenstein dimension (see \cite[Corollaries 4.4.6 and 4.4.13]{GorensteinDimensions}). 

\begin{theorem} \label{seqexataelementsoemAuseBassClass} \cite[Corollary 2.9]{HomologicalAspectsOfSemidualizingModules} The class $\mathcal{A}_C(R)$ 
        contains every $R$-module of finite $C$-injective dimension. 
\end{theorem}
\begin{theorem}\cite[Theorem 2.11]{HomologicalAspectsOfSemidualizingModules} \label{formulaIcPc} Let $M$ be an $R$-module. Then $\mathcal{I}_C\operatorname{-id}_R(M)=\operatorname{id}_R(C \otimes_R M)$.
\end{theorem}

We refer the reader to \cite{SemidualizingModules, HomologicalAspectsOfSemidualizingModules} for detailed results about semidualizing module and $C$-injective dimension.

If $R$ is a Cohen-Macaulay ring with a canonical module $\omega_R$, we can define $\omega_C$ as $\Hom_R(C, \omega_R)$. Clearly, when $C=R$ we observe that $\omega_C$ coincides with $\omega_R$. The module $\omega_C$ has been studied by Bagheri and Taherizadeh in \cite{CCanonicalModule}. 

\begin{proposition}
Let $R$ be a Cohen-Macaulay local ring with a canonical module $\omega_R$. Then $\omega_C$ is a semidualizing $R$-module. In particular, $\omega_C$ is indecomposable.
\end{proposition}
\begin{proof}
    This follows from \cite[Remark 4.4]{CCanonicalModule} and Proposition \ref{Propbasicasemidualizantes}(2).
\end{proof}

\begin{theorem}\cite[Theorem 4.9]{CCanonicalModule}
\label{copiaswc}
    Let $R$ be a Cohen-Macaulay local ring with a canonical module $\omega_R$. If $M$ is a maximal Cohen-Macaulay $R$-module with  $\mathcal{I}_C\operatorname{-id}_R(M)<\infty$, then $M\cong \omega_C^t$ for some positive number $t$.
\end{theorem}

\section{Extension of results to $C$-injective dimension of $\Hom$}

In \cite{injectivedimensiontakahashi}, Ghosh and Takahashi proved the following: 
Let $M$ be an $R$-module, and $N$ be a maximal Cohen-Macaulay $R$-module. Suppose that $\Hom_R(M,N)$ has finite injective dimension. 
\begin{enumerate}
\item If $\Hom_R(M,N)\not=0$ and $\Ext_R^i(M,N)=0$ for all $1\leq i \leq d-1$, then  $R$ is Cohen-Macaulay,  $N$ has finite injective dimension and $M\cong \Gamma_\mathfrak{m}(M)\oplus R^r$ for some $r\geq 0$.
\item If $M$ and $N$ are nonzero and  $\Ext_R^i(M,N)=0$ for all $1\leq i \leq d$, then $M$ is free and $N$ has finite injective dimension.
\end{enumerate}

In this section, our objective is to extend these results and some results from Section 3 regarding the $C$-injective dimension of $\Hom$. To achieve these goals, we give the following lemma that will be the key for the proofs.

\begin{lemma}\label{Lemainteresante}
    Let $M$ and $N$ be two $R$-modules. Suppose $t$ is a nonnegative integer. If $\Hom_R(M,N)\in \mathcal{A}_C(R)$ and $N \in \mathcal{A}_C(R)$, then:
    \begin{enumerate}
        \item[(1)] The natural map $$\Phi_M: C \otimes_R \Hom_R(M,N)\to \Hom_R(M, C \otimes_R N), c \otimes \varphi \mapsto \varphi_c,$$
        where $\varphi_c: M \to C \otimes_R N$ is given by $m \mapsto c \otimes \varphi(m)$, is an isomorphism.
        \item[(2)]  If $\Ext_R^i(M,N)=0$ for all $1\leq i \leq t$, then $\Ext_R^i(M, C \otimes_R N)=0$ for all $1\leq i \leq t.$
    \end{enumerate}
\end{lemma}
\begin{proof}
$(1)$ If $M=R^r$ for some $r\geq 0$, it is clear that $\Phi_M$ is an isomorphism. For the general case, we consider a finite presentation of $M$,
$$R^s \stackrel{\partial_1}{\longrightarrow} R^r \stackrel{\partial_0}{\longrightarrow} M \longrightarrow 0 .$$
Applying $\Hom_R(-,N)$ to this sequence, we obtain an exact sequence
$$\xymatrix{0 \ar[r]& \Hom_R(M,N) \ar[r]^{\partial_0^\ast}& \Hom_R(R^{r}, N)  \ar[r]^{\partial_1^\ast}& \Hom_R(R^{s}, N)}.$$
Since $\operatorname{Hom}_R(M,N)$ and $N$ belong to $\mathcal{A}_C(R)$, all the $R$-modules in the exact sequence above also belong to $\mathcal{A}_C(R)$.

Thus, by applying \cite[Lemma 3.1.12]{SemidualizingModules}, we have the following commutative diagram,
$$\xymatrix{0 \ar[r]& C \otimes_R \Hom_R(M,N) \ar[d]_{\Phi_M} \ar[r]^{C \otimes_R \partial_0^\ast}& C \otimes_R \Hom_R(R^{r}, N)  \ar[d]^{\Phi_{R^{r}}} \ar[r]^{C\otimes_R \partial_1^\ast}& C \otimes_R \Hom_R(R^{s}, N) \ar[d]^{\Phi_{R^{s}}} \\
0 \ar[r]& \Hom_R(M, C \otimes_R N) \ar[r]& \Hom_R(R^{r}, C \otimes_R N) \ar[r]& \Hom_R(R^{s}, C \otimes_R N)
},$$
where the rows are exact. Since $\Phi_{R^{r}}$ and $\Phi_{R^{s}}$ are isomorphisms, it follows from the Five Lemma that $\Phi_M$ is an isomorphism.

$(2)$ Let 
$$\cdots \longrightarrow R^{n_{i+1}} \longrightarrow R^{n_i} \longrightarrow \cdots \longrightarrow R^{n_1} \longrightarrow R^{n_0} \longrightarrow M \longrightarrow 0 $$
be a minimal free resolution of $M$.
Since $\Ext_R^i(M,N)=0$ for all $1\leq i \leq t$, we have the following exact sequence:
$$0 \longrightarrow \operatorname{Hom}_R(M, N) \longrightarrow \operatorname{Hom}_R\left(R^{n_0}, N\right) \longrightarrow \cdots \longrightarrow \operatorname{Hom}_R\left(R^{n_{t+1}}, N\right).$$
As $\Hom_R(M,N)$ and $N$ belong to $\mathcal{A}_C(R)$, all the $R$-modules in the exact sequence above  also belong to $\mathcal{A}_C(R)$. Thus, by using (1) and \cite[Lemma 3.1.12]{SemidualizingModules}, we obtain the following commutative diagram:
$$\tiny \xymatrix{0 \ar[r]& C \otimes_R \Hom_R(M,N) \ar[r] \ar[d]_{\Phi_M}& C \otimes_R \Hom_R(R^{n_0},N) \ar[r] \ar[d]^{\Phi_{R^{n_0}}}& \cdots \ar[r]&C \otimes_R \Hom_R(R^{n_{t+1}}, N) \ar[d]^{\Phi_{R^{n_{t+1}}}} \\ 0 \ar[r]& \Hom_R(M, C \otimes_R N) \ar[r]& \Hom_R(R^{n_0}, C \otimes_R N) \ar[r]&  \cdots \ar[r]& \Hom_R(R^{n_{t+1}},C \otimes_R N)   
},$$
where the first row is exact and the vertical maps are isomorphisms. Hence, the second row is also exact and we conclude that $\Ext_R^i(M,C \otimes_R N)=0$ for  all $1\leq i \leq t$. 
\end{proof}

\begin{remark}\label{lemaparaaplicar}
    Let $M$ and $N$ be $R$-modules. 
    \begin{enumerate}
        \item[(1)] If $N\not=0$, then $C \otimes_R N\not=0$. 
        \item[(2)] If $N \in \mathcal{A}_C(R)$, then $\depth(N)=\depth(C \otimes_R N)$. In particular, if $N \in \mathcal{A}_C(R)$ is maximal Cohen-Macaulay, then $C \otimes_R N$ is maximal Cohen-Macaulay.
        \item[(3)] If $N \in \mathcal{A}_C(R)$, $\Hom_R(M,N)$ is an (resp. a nonzero) $R$-module of finite  $C$-injective dimension and $\Ext_R^i(M,N)=0$ for all $1\leq i \leq t$. Then  $\Hom_R(M,C\otimes_R N)$ is an (resp. a nonzero) $R$-module of finite injective dimension and  $\Ext_R^i(M, C \otimes_R N)=0$ for all $1\leq i \leq t.$
    \end{enumerate}
\end{remark}
\begin{proof}
Part $(1)$ is trivial, and part $(2)$ follows from \cite[Lemma 2.11(1)]{LinkageOfModulesAndTheSerreConditions}.

$(3)$ Since $\Hom_R(M,N)$ has finite $C$-injective dimension, we can apply Theorem \ref{seqexataelementsoemAuseBassClass} to conclude that $\Hom_R(M,N) \in \mathcal{A}_C(R)$. Additionally, by Theorem \ref{formulaIcPc}, we have $\operatorname{id}_R(C \otimes_R \Hom_R(M,N))<\infty$.

Applying Lemma \ref{Lemainteresante}, we obtain that 
\begin{equation}\label{isoChOM}
C \otimes_R \Hom_R(M,N) \cong \Hom_R(M, C \otimes_R N) 
\end{equation}
and $\Ext_R^i(M, C\otimes_R N)=0 \mbox{ for all } 1\leq i \leq t.$

Since $\operatorname{id}_R(C \otimes_R \Hom_R(M,N))<\infty$, the isomorphism \eqref{isoChOM} implies that  $\operatorname{id}_R(\Hom_R(M, C \otimes_R N))<\infty$.

Additionally, if $\Hom_R(M,N)\neq 0$, then from \eqref{isoChOM} we see that $\Hom_R(M, C \otimes_R N)\neq 0$.
\end{proof}

The following result provides an extension of \cite[Theorem 2.3]{injectivedimensiontakahashi}.

\begin{theorem}\label{generalizacionteo2.3}
Let  $M$ and $N$ be $R$-modules such that $N \in \mathcal{A}_C(R)$ is  maximal Cohen-Macaulay. Suppose that $\Hom_R(M,N)$ is nonzero and has finite $C$-injective dimension and that $\Ext_R^i(M,N)=0$ for all $1\leq i \leq d-1$. Then $R$ is a Cohen-Macaulay ring, $N$ has finite $C$-injective dimension and $M \cong \Gamma_\mathfrak{m}(M) \oplus R^r$ for some  $r\geq 0$.
\end{theorem}
\begin{proof}
By Remark \ref{lemaparaaplicar}, we have that $C \otimes_R N$ is  maximal Cohen-Macaulay,  $\Hom_R(M,C \otimes_R N)$ is a nonzero $R$-module of finite injective dimension, and  $\Ext_R^i(M,C \otimes_R N)=0$ for all $1\leq i \leq d-1$. Using \cite[Theorem 2.3]{injectivedimensiontakahashi}, we deduce that $R$ is a Cohen-Macaulay ring, $M\cong \Gamma_\mathfrak{m}(M) \oplus R^r$ for some nonnegative integer $r$, and $C \otimes_R N$ has finite injective dimension. Furthermore, Theorem \ref{formulaIcPc} implies that $N$ has finite $C$-injective dimension.
\end{proof}
The following corollary generalizes \cite[Corollary 2.10(1)]{injectivedimensiontakahashi}.
\begin{corollary} \label{generalizaciondecor2.10(1)}
 Let $M$ be a nonzero $R$-module such that $\Ext_R^i(M,R)=0$ for all $1\leq i \leq d-1$,  and $M^\ast$ a nonzero $R$-module of finite $C$-injective dimension.
Then $R$ is Cohen-Macaulay, $C$ is a canonical module of $R$ and $M \cong \Gamma_\mathfrak{m}(M) \oplus R^r$ for some $r\geq 0$.
\end{corollary}
\begin{proof}
    By taking $N=R$ in Theorem \ref{generalizacionteo2.3}, it follows that $R$ is Cohen-Macaulay, $M \cong \Gamma_\mathfrak{m}(M)\oplus R^r$ for some $r\geq 0$, and $\mathcal{I}_C\operatorname{-id}_R(R)<\infty$. According to Theorem \ref{formulaIcPc}, $\operatorname{id}_R(C)<\infty$, which implies that $C$ is a canonical module of $R$.
\end{proof}

The following theorem extends \cite[Theorem 2.5]{injectivedimensiontakahashi}.
 \begin{theorem}\label{generalizaciondoteo2.5}
     Let $M$ and $N$ be nonzero $R$-modules such that $N \in \mathcal{A}_C(R)$ is maximal Cohen-Macaulay. Suppose that $\Hom_R(M,N)$ has finite $C$-injective dimension and that $\Ext_R^i(M,N)=0$ for all $1\leq i \leq d$. Then $R$ is Cohen-Macaulay, $M$ is free and $N$ has finite $C$-injective dimension. 
 \end{theorem}
\begin{proof} 
By Remark \ref{lemaparaaplicar}, we know that $C \otimes_R N$ is maximal Cohen-Macaulay, $\Hom_R(M,C \otimes_R N)$ has finite injective dimension, and $\Ext_R^i(M,C \otimes_R N)=0$ for all $1 \leq i \leq d-1$. Now, by \cite[Theorem 2.5]{injectivedimensiontakahashi}, $C \otimes_R N$ has finite injective dimension, and $M$ is free. Theorem \ref{formulaIcPc} implies that $N$ has finite $C$-injective dimension, and Bass' theorem implies that $R$ is Cohen-Macaulay.
\end{proof}

Next, we provide one of the main results of this section that generalizes \cite[Corollary 2.10(2)]{injectivedimensiontakahashi}.

\begin{corollary}\label{generalizacionCorolario2.10(2)}
    Let  $M$ be a nonzero $R$-module. If $M^\ast$ has finite $C$-injective dimension and $\Ext_R^i(M,R)=0$ for all $1\leq i \leq d$, then $R$ is Cohen-Macaulay, $C$ is a canonical module of $R$ and $M$ is free.
\end{corollary}
\begin{proof}
Take $N=R$ in Theorem \ref{generalizaciondoteo2.5}.
\end{proof}

Next, by Corollary \ref{generalizacionCorolario2.10(2)}, we provide an affirmative answer to the Auslander-Reiten conjecture when $\operatorname{Hom}_R(M,R)$ has finite $C$-injective dimension.

\begin{corollary} \label{ARCparaCidimHom(M,R)} The Auslander-Reiten conjecture holds true for (finitely generated) $R$-modules $M$ such that $\mathcal{I}_C\operatorname{-id}_R(\Hom_R(M,R))<\infty$.
\end{corollary}

As a consequence of Theorem \ref{generalizaciondoteo2.5}, we have the following interesting corollaries that determine when a module is free.

\begin{corollary}\label{corpdGprojective}
    Let $R$ be a Cohen-Macaulay local ring, and $M$ and $N$ be nonzero $R$-modules such that $N$ is $G$-projective. Suppose that $\operatorname{pd}_R(\Hom_R(M,N))<\infty$ and that  $\Ext_R^i(M,N)=0$ for all $1\leq i \leq d$. Then $M$ and $N$ are free.
\end{corollary}
\begin{proof} 
We may assume that $R$ is complete, and therefore, $R$ admits a canonical module $\omega_R$. Since $N$ is $G$-projective, then $N \in \mathcal{A}_{\omega_R}(R)$, and by the Auslander-Bridger formula, $N$ is maximal Cohen-Macaulay. Now, by \cite[Theorem 3.4.11]{GorensteinDimensions}, we obtain that $\mathcal{I}_{\omega_R}\operatorname{-id}_R(\operatorname{Hom}_R(M,N))= \operatorname{id}_R( \omega_R \otimes_R \operatorname{Hom}_R(M,N))<\infty$. It follows from Theorem \ref{generalizaciondoteo2.5} that $M$ is free. Thus, $\operatorname{pd}_R(N)<\infty$ since $\operatorname{pd}_R(\Hom_R(M,N))<\infty$, and by the Auslander-Buchsbaum formula, $N$ is free.
\end{proof}

\begin{corollary}\label{ARCpdHom(M,R)}
    Let $R$ be a Cohen-Macaulay local ring. Then the Auslander-Reiten conjecture holds true for (finitely generated) $R$-modules $M$ such that $\operatorname{pd}_R(\Hom_R(M,R))<\infty$.
\end{corollary}

In the following example, we can see that the assumption $N \in \mathcal{A}_C(R)$ in Theorem \ref{generalizaciondoteo2.5} cannot be removed.

\begin{example} Let $R$ be a Cohen-Macaulay local ring that is not Gorenstein, with a canonical module $\omega_R$. Note that $\Hom_R(\omega_R, \omega_R)\cong R$ has finite $\omega_R$-injective dimension since $\operatorname{id}_R(\omega_R)<\infty$. It is clear that $\Ext_R^i(\omega_R, \omega_R)=0$ for all $i\geq 1$. However, it is important to observe that $\omega_R$ is not a free module because $R$ is not Gorenstein.
\end{example}

Using Remark \ref{lemaparaaplicar}, Theorem \ref{otrageneralizacionde2.15} and Theorem \ref{formulaIcPc}, we get a theorem with the same implication of Theorem \ref{otrageneralizacionde2.15}, but for a larger class of rings and modules.

\begin{theorem}\label{CorolariogeneralizadoparaIc}
Let $M$ and $N$ be nonzero $R$-modules such that $N \in \mathcal{A}_C(R)$, and let $t=\depth(N)$. Suppose that $\Hom_R(M,N)$ has finite $C$-injective dimension and that $\Ext_R^i(M,N)=\Ext_R^j(M,R)=0$  for all $1\leq i \leq t$ and  $1\leq j \leq d$. Then $R$ is Cohen-Macaulay, $M$ is free and  $N$ has finite $C$-injective dimension. 
\end{theorem}
\begin{proof} 
By Remark \ref{lemaparaaplicar}, $C\otimes_R N$ is nonzero, $\depth(C\otimes_R N)=t$,     $\Hom_R(M, C\otimes_R N)$ has finite injective dimension and $\Ext_R^i(M,C \otimes_R N)=0$ for all $1\leq i \leq t$. Thus, by  Theorem \ref{otrageneralizacionde2.15}, $M$ is free and $C \otimes_R N$ has finite injective dimension. Theorem \ref{formulaIcPc} implies that $N$ has finite $C$-injective dimension, and Bass' theorem implies that $R$ is Cohen-Macaulay.
\end{proof}
\begin{corollary}\label{otrocritpd}
Let $R$ be a Cohen-Macaulay local ring, and let $M$ and $N$ be  nonzero $R$-modules such that  $\operatorname{G-dim}_R(N)<\infty$. Let $t=\depth(N)$. Suppose that  $\operatorname{pd}_R(\Hom_R(M,N))<\infty$, and that $\Ext_R^i(M,N)=\Ext_R^j(M,R)=0$ for all $1\leq i \leq t$ and  $1\leq j\leq d$.   Then $M$ is free and $\operatorname{pd}_R(N)<\infty.$
\end{corollary}

\begin{proof}
We may assume that $R$ is complete, and hence $R$ admits a canonical module $\omega_R$. As $\operatorname{G-dim}_R(N)<\infty$, then $N \in \mathcal{A}_{\omega_R}(R)$.   On the other hand, by \cite[Theorem 3.4.11]{GorensteinDimensions},  we have that $\mathcal{I}_{\omega_R}\operatorname{-id}_R(\operatorname{Hom}_R(M,N))=\operatorname{id}_R( \omega_R \otimes_R \operatorname{Hom}_R(M,N))<\infty
$. It follows from Theorem \ref{CorolariogeneralizadoparaIc} that $M$ is free. Since $\Hom_R(M,N)$ has finite projective, $N$ does as well.
\end{proof}

Similarly to the proof of Theorem \ref{estasigeneralizacionde215}, we have the following result that extends such theorem.
\begin{theorem}\label{2.15generalizadoIC}
Let $M$ be a nonzero $R$-module such that $M \in \mathcal{A}_C(R)$.  Suppose that
 $\Hom_R(M,M)$ has finite $C$-injective dimension, and that $\Ext_R^i(M,M)=\Ext_R^j(M,R)=0$ for all $1\leq i \leq d-1$ and $1\leq j \leq d$.  Then  $R$ is Cohen-Macaulay,  $M$ is free and $C$ is a canonical module for $R$. 
\end{theorem}
\begin{proof} 
We may assume that $R$ is complete. As $\Hom_R(M,M)$ is a nonzero $R$-module of finite $C$-injective dimension, by Theorem \ref{formulaIcPc} and  Bass' theorem, $R$ is Cohen-Macaulay, and hence admits a canonical module $\omega_R$.

Let $t=\depth(M)$. We claim that $M$ is maximal Cohen-Macaulay, i.e., $t=d$. Suppose on the contrary that $t\leq d-1$. Since the vanishing assumption shows $\Ext_R^i(M,M)=0$ for all $1\leq i \leq t$, according to Theorem  \ref{CorolariogeneralizadoparaIc}, $M$ is free, which contradicts $t<d$. Hence, we conclude that $t=d$.

By Theorem \ref{generalizacionteo2.3}, we have $M \cong \Gamma_{\mathfrak{m}}(M) \oplus R^r$ for some $r\geq 0$ and $\mathcal{I}_C\operatorname{-id}_R(M)<\infty$. However, since $M$ is maximal Cohen-Macaulay, we have $\Ext_R^d(M, \omega_R)=0$. Applying the local duality theorem, we find that $\Gamma_{\mathfrak{m}}(M)=0$, which implies $M \cong R^r$.

Since $\mathcal{I}_C\operatorname{-id}_R(M)<\infty$ and $M$ is free, by Theorem \ref{formulaIcPc}, we conclude that $\operatorname{id}_R(C)<\infty$. Hence, $C$ is a canonical module for $R$. Therefore, the theorem is proved.

\end{proof}

\begin{corollary}\label{ARCforCidHom(M,M)}
  The Auslander-Reiten conjecture holds true for (finitely generated) $R$-modules $M \in \mathcal{A}_C(R)$ with $\mathcal{I}_C\operatorname{-id}_R( \Hom_R(M,M))<\infty$.
\end{corollary}

As an application of Theorem \ref{2.15generalizadoIC}, we derive the following theorem which generalizes \cite[Corollary 2.14]{injectivedimensiontakahashi} and give a criteria for a module of finite Gorenstein dimension to be free over a Cohen-Macaulay ring.

\begin{theorem}\label{pdHom}
    Let $R$ be a Cohen-Macaulay local ring and $M$ be an $R$-module with \linebreak
    $\operatorname{G-dim}_R(M)<\infty$. The following are equivalent:
\begin{enumerate}
    \item $M$ is free.
    \item $\Hom_R(M,M)$ is free and $\Ext_R^i(M,M)=0$ for all $1\leq i \leq d$.
    \item $\Hom_R(M,M)$ has finite projective dimension and $\Ext_R^i(M,M)= \Ext_R^j(M,R) =0$ for all $1\leq i \leq d-1$ and $1\leq j \leq d$.
\end{enumerate}
\end{theorem}
\begin{proof}
The implication $(1) \Rightarrow (2)$ is trivial.

$(2) \Rightarrow (3)$. We only need to show that $\Ext_R^j(M,R)=0$ for all $1\leq j \leq d$. Let us assume that $M\neq 0$. Since $\Hom_R(M,M)$ is free, we get that $\Hom_R(M,M)$ is maximal Cohen-Macaulay. Thus, as $\Ext_R^i(M,M)=0$ for all $1\leq i \leq d$, we conclude from Lemma \ref{dephdeHom}(2) that $M$ is maximal Cohen-Macaulay. Applying the Auslander-Bridger formula, we deduce that $M$ is $G$-projective, and therefore $\Ext_R^j(M,R)=0$ for all $j>0$.

$(3) \Rightarrow (1)$. We may assume that $R$ is complete, and hence $R$ admits a canonical module $\omega_R$.  As $\operatorname{G-dim}_R(M)<\infty$, then $M \in \mathcal{A}_{\omega_R}(R)$. On the other hand, since $\operatorname{pd}_R(\Hom_R(M,M))<\infty$, it follows from \cite[Theorem 3.4.11]{GorensteinDimensions} that $\Hom_R(M,M)$ has finite $\omega_R$-injective dimension.  Hence, by Theorem \ref{2.15generalizadoIC}, $M$ is free. 
\end{proof}

Next, according to Theorem \ref{pdHom}, in the case where $R$ is Cohen-Macaulay, we can provide a positive answer to the Auslander-Reiten conjecture when $M$ has finite Gorenstein dimension and $\operatorname{Hom}_R(M,M)$ has finite projective dimension.
\begin{corollary}\label{ARGorenstein}
   Let $R$ be a Cohen-Macaulay local ring. Then the Auslander-Reiten conjecture holds true for (finitely generated) $R$-modules $M$  with  $\operatorname{G-dim}_R(M)<\infty$ and  $\operatorname{pd}_R(\Hom_R(M,M))<\infty$.
\end{corollary}

\begin{remark}\label{rmaa}
    Corollaries \ref{corpdGprojective}, \ref{ARCpdHom(M,R)},  \ref{ARGorenstein} and Theorem \ref{pdHom} were obtained as consequences of the some theorems of this section. However, these results can be derived from \cite{HomandExtRevisited} with significantly weaker hypotheses, and in some of them relaxing the hypotheses of finite $\operatorname{G}$-dimension. For instance, Corollary \ref{corpdGprojective} (from Corollary \ref{ARCpdHom(M,R)} follows directly) can be obtained from \cite[Lemmas 3.1 and 3.3]{HomandExtRevisited}, noting that the assumptions immediately force $\Hom_R(M, N)$ to be maximal Cohen-Macaulay by the depth lemma or Lemma \ref{dephdeHom}. Similarly, for Theorem \ref{pdHom} (from which Corollary \ref{ARGorenstein} follows directly), we may replace the assumption that $M$ has finite $G$-dimension with the weaker assumption that $\Omega_R^{d-\text{depth}_R(M)}(M)$ is a syzygy of a maximal Cohen-Macaulay $R$-module and that given  $d-\operatorname{depth}_R(M)\leq i \leq d$, the nonvanishing of $\text{Ext}^i_R(M, R)$ occurs if and only if $i = d - \text{depth}_R(M)$. With this hypothesis instead, (1) $\Rightarrow$ (2) and (2) $\Rightarrow$ (3) follow the same with the key point being the standard argument that $M$ is maximal Cohen-Macaulay under the hypotheses of (2). (3) $\Rightarrow$ (1) follows since the assumption $\text{Ext}^i_R(M, R) = 0$ for all $1 \leq i \leq d$ then forces $M$ to be maximal Cohen-Macaulay in consideration of above. The assumption $\text{Ext}^i_R(M, M) = 0$ for all $1 \leq i \leq d - 1$ then forces $\Hom_R(M, M)$ to be maximal Cohen-Macaulay by the depth lemma whence it is free, and we may then apply \cite[Theorem 3.8]{HomandExtRevisited} to conclude that $M$ is free. The same weakening of the finite G-dimension hypothesis can be applied to Corollary \ref{ARGorenstein}, but in either case Corollary
\ref{ARGorenstein} is an immediate consequence of \cite[Theorem 3.8]{HomandExtRevisited}. On the other hand, in an independent work, Ghosh and Dey \cite{FiniteHomologicalDimensionOfHomv2} showed that in these results the assumption of Cohen-Macaulayness on $R$ can be removed and that, additionally, in Corollary \ref{corpdGprojective} and Theorem \ref{pdHom}  the number of the vanishing can be reduced. In view of all this, these results recover some of \cite{HomandExtRevisited, FiniteHomologicalDimensionOfHomv2} under slightly more stringent hypotheses.
\end{remark} 

\section{Characterizations for a semidualizing module to be canonical module}
In this section, we give some results that establish criteria for a semidualizing $R$-module $C$ to be a canonical module for $R$ in terms of finite $C$-injective dimension of certain Hom. Furthermore, we generalize and improve the results given by Ghosh and Takahashi in \cite[Proposition 3.2, Theorem 3.6]{injectivedimensiontakahashi}.

\begin{proposition} \label{generalizaciondeProposicion32} Let $M$ be an $R$-module.
If  $\depth( \Hom_R(M,M))=d$ and $\Hom_R(M,M)$ has finite $C$-injective dimension, then $R$ is Cohen-Macaulay and $C$ is a canonical module for $R$.
\end{proposition}
\begin{proof}
  We can assume that $R$ is complete. The equality $\depth(\Hom_R(M,M))=d$ shows that $\Hom_R(M,M)\neq 0$. Since $\mathcal{I}_C\operatorname{-id}_R(\Hom_R(M,M))<\infty$, according to Theorem \ref{formulaIcPc}, $C \otimes_R \Hom_R(M,M)$ is a nonzero $R$-module of finite injective dimension. Therefore, by Bass' theorem, it follows that $R$ is Cohen-Macaulay, and hence $R$ admits a canonical module $\omega_R$. Moreover, Theorem \ref{copiaswc} implies that $\Hom_R(M,M)\cong \omega_C^n$ for some $n\geq 1$. Then
    \begin{align*}
        R^{n^2}&\cong \Hom_R(\omega_C, \omega_C)^{n^2} \\
        & \cong \Hom_R( \omega_C^n, \omega_C^n)\\ & \cong \Hom_R( \Hom_R(M,M), \Hom_R(M,M)) \\
        & \cong \Hom_R(M \otimes \Hom_R(M,M), M),
    \end{align*}
    and hence $ R^{n^2} \cong \Hom_R(M \otimes \Hom_R(M,M), M).$ Now, consider the map  $f: M \to M \otimes_R \Hom_R(M,M)$ defined by $f(m)=m \otimes \operatorname{Id}_M$. 
This map is a split injection (see \cite[Lemma 3.1]{MCMtensorproductsandvanishingofExtmodules}), so there exists an $R$-module $N$ such that 
$M \oplus N \cong M \otimes_R \Hom_R(M,M).$ Then
\begin{align*}
\omega_C^n \oplus \Hom_R(N,M) & \cong \Hom_R(M,M) \oplus \Hom_R(N,M)\\ &\cong \Hom_R(M \oplus N, M) \\ & \cong \Hom_R ( M \otimes_R \Hom_R(M,M),M)\\ & \cong R^{n^2}.
\end{align*}  
Thus, $\omega_C^n \oplus \Hom_R(M,N)\cong R^{n^2}$. 
As $R$ is indecomposable and satisfies the Krull-Smichdt Theorem (\cite[Proposition 1.18]{livroMCMyoshino}), then  $\omega_C \cong R.$ Hence,
\begin{align*}
    C &\cong \Hom_R(\Hom_R(C,\omega_R), \omega_R) \\& \cong \Hom_R(\omega_C,\omega_R)\\ & \cong   \Hom_R(R,\omega_R)\\& \cong  \omega_R.
\end{align*}
\end{proof}

\begin{corollary}\label{corolarin}
 Let  $M$ be a  maximal Cohen-Macaulay $R$-module such that  $\Hom_R(M,M)$ has finite $C$-injective dimension. If $\Ext_R^i(M,M)=0$
     for all $1\leq i \leq d-1$, then 
      $R$ is Cohen-Macaulay  and $C$ is a canonical module for $R$.
\end{corollary}
\begin{proof} As  $M$ is maximal Cohen-Macaulay, then  $M\not=0$, and hence $\Hom_R(M,M)\not=0$. Since \linebreak  $\mathcal{I}_C\operatorname{-id}_R(\Hom_R(M,M))<\infty$, according to Theorem \ref{formulaIcPc}, $C \otimes_R \Hom_R(M,M)$ is a nonzero $R$-module of finite injective dimension. Therefore, by Bass' theorem, it follows that $R$ is Cohen-Macaulay. Now, consider a minimal free resolution of $M$ given by
$$\cdots \to R^{n_3} \to R^{n_2} \to R^{n_0} \to M \to 0.$$
Since $\Ext_R^i(M,M)=0$ for all $1\leq i \leq d-1$, we obtain an exact sequence
$$0 \to \Hom_R(M,M) \to M^{n_0} \to \cdots \to M^{n_e},$$
where $e=\max\{1,d\}$. According to the depth lemma, $\Hom_R(M,M)$ is  maximal Cohen-Macaulay. Consequently, by Proposition \ref{generalizaciondeProposicion32}, $C$ is a canonical module for $R$.
\end{proof}

The following generalizes \cite[Theorem 3.6]{injectivedimensiontakahashi}.
\begin{theorem}\label{caracterizacionessemidualizingserdualizing} The following statements are equivalent:
\begin{enumerate}
    \item[(1)] $R$ is Cohen-Macaulay and 
    $C$ is a canonical module for $R$.
    \item[(2)] $R$ admits a module $M$ such that $\Ext_R^{j}(M,R)=0$ for all $1\leq j \leq d-1$ and $M^\ast$ is nonzero of finite $C$-injective dimension.
    \item[(3)] $R$ admits a nonzero module  $M \in \mathcal{A}_C(R)$ such that $\mathcal{I}_C\operatorname{-id}_R( \Hom_R(M,M))<\infty$, and  \linebreak $\Ext_R ^i(M,M)=\Ext_R^{j}(M,R)=0$ for all $1\leq i \leq d-1$ and $1\leq j \leq d$.
    \item[(4)]  $R$ admits a maximal Cohen-Macaulay module  $M$ such that $\Ext_R^i(M,M)=0$ for all $1\leq i \leq d-1$ and $\mathcal{I}_C\operatorname{-id}_R(\Hom_R(M,M))<\infty$.
    \item[(5)] $R$ admits a  module $M$ such that $\Hom_R(M,M)$ is maximal Cohen-Macaulay and \linebreak  $\mathcal{I}_C\operatorname{-id}_R(\Hom_R(M,M))<\infty.$
\end{enumerate}  
\end{theorem}
\begin{proof}
Note that the implications $(1)\Rightarrow (2),(3),(4),(5)$ hold when we take $M=R$. The reverse implications $(2)\Rightarrow(1)$, $(3)\Rightarrow (1)$, $(4)\Rightarrow (1)$, and $(5) \Rightarrow (1)$ follow from Corollary \ref{generalizaciondecor2.10(1)}, Theorem \ref{2.15generalizadoIC}, Corollary \ref{corolarin}, and Proposition \ref{generalizaciondeProposicion32}, respectively.

\end{proof}

By setting $C=R$ in Theorem \ref{caracterizacionessemidualizingserdualizing},  we obtain the following result that improves \cite[Theorem 3.6]{injectivedimensiontakahashi} and provides characterizations for a Gorenstein local ring.
\begin{corollary}\label{Gorentein}
The following statements are equivalent:
\begin{enumerate}
   \item[(1)] $R$ is Gorenstein. 
   \item[(2)] $R$ admits a module $M$ such that $\Ext_R^{j}(M,R)=0$ for all $1\leq j \leq d-1$ and $M^\ast$ is nonzero of finite injective dimension.
    \item[(3)] $R$ admits a nonzero module  $M$ such that $\Ext_R^i(M,M)=\Ext_R^{j}(M,R)=0$ for all  $1\leq i \leq d-1$ and $1\leq j \leq d$, and $\operatorname{id}_R(\Hom_R(M,M))<\infty$.
    \item[(4)] $R$ admits a  maximal Cohen-Macaulay module $M$ such that $\Ext_R^i(M,M)=0$ for all $1\leq i \leq d-1$ and $\operatorname{id}_R(\Hom_R(M,M))<\infty$.
    \item[(5)] $R$ admits a module $M$ such that $\Hom_R(M,M)$ is maximal Cohen-Macaulay and \linebreak  $\operatorname{id}_R(\Hom_R(M,M))<\infty$.
\end{enumerate}  
\end{corollary}
\section*{Acknowledgment}
The authors acknowledge the anonymous referee for reading the paper carefully and giving valuable comments. In particular, we thank him/her for indicating how some results presented in this paper recover some results of \cite{HomandExtRevisited} under slightly more hypotheses, which gave rise to Remark \ref{rmaa}.


\begin{thebibliography}{999999}
\bibitem{AnneauxdeGorensteinettorsionenalgebrecommutative} { M. Auslander}, \textsl{Anneaux de Gorenstein et torsion en algèbre commutative}, Séminaire d'algèbre commutative 1966/67, notes by {  M. Mangeney, C. Peskine  L. Szpiro}, École Normale Supérieure de Jeunes Filles, Paris, 1967.
\bibitem{Stablemoduletheory} { M. Auslander, M. Bridger}, \textsl{Stable module theory}, Mem. of the AMS 94, Amer. Math. Soc., Providence, R.I.,1969.
\bibitem{OnageneralizedversionoftheNakayamaconjecture} { M. Auslander, I. Reiten}, On a generalized version of the Nakayama conjecture, Proc. Amer. Math. Soc. \textbf{52} (1975)  69--74.
\bibitem{CCanonicalModule} { M. Bagheri, A. Taherizadeh}, $C$\--canonical modules, 
§Int. Electron. J. Algebra  \textbf{30} (2021), 243--259. 
\bibitem{bruns} {W. Bruns, J. Herzog}, \textsl{Cohen-Macaulay rings}, Sec. Edt., Cambridge Univ. Press, 1998.
\bibitem{GorensteinDimensions} { L. W. Christensen}, \textsl{Gorenstein Dimensions}, Lect. Notes Math. \textbf{1747}, Springer-Verlag, Berlin, 2000.

\bibitem{HomandExtRevisited}{  H. Dao, M.  Eghbali, J. Lyle}, Hom and Ext, revisited, Journal of Algebra \textbf{571} (2021), 75--93.
\bibitem{SomeHomolocalPropertiesofModulesoveraCompleteIntersectionwithApplications} { H. Dao}, Some Homological Properties of Modules over a Complete Intersection, with Applications, in 
 \textsl{Commutative Algebra: Expository Papers Dedicated To David Eisenbud On The Occasion Of His 65th Birthday}, 335--371, Springer, New York, 
 2013.
\bibitem{FiniteHomologicalDimensionOfHomv2}
S. Dey and D. Ghosh, Finite homological dimension of {H}om and vanishing of {E}xt, arXiv preprint arXiv:2310.10607v1, 2023.

\bibitem{LinkageOfModulesAndTheSerreConditions} { M. Dibaei, A. Sadeghi},  Linkage of modules and the Serre conditions, J. Pure Appl. Algebra \textbf{219} (2015), 4458--4478.
\bibitem{GorensteinModulesAndRelatedModules} { H. Foxby},  Gorenstein modules and related modules, Math. Scand. \textbf{31} (1972), 267--284. 

\bibitem{injectivedimensiontakahashi} { D. Ghosh, R. Takahashi}, Auslander-Reiten conjecture and finite injective dimension of Hom, Kyoto J. Math. \textbf{64} (2024), 229--243.

\bibitem{GdimensionandGeneralizedPerfectIdeals} { E. Golod},  G-dimension and generalized perfect ideals,  Trudy Mat. Inst. Steklov \textbf{165} (1984),  62--66.
\bibitem{AuslanderReitenConjectureforNormalRings}
{K. Kimura}, Auslander-{R}eiten conjecture for normal rings, arXiv preprint arXiv:2304.03956v1, 2023.
\bibitem{MCMtensorproductsandvanishingofExtmodules}{  K. Kimura, Y. Otake, R. Takahashi}, Maximal Cohen-Macaulay tensor products and vanishing of Ext modules, Bull. Lond. Math. Soc., \textbf{54} (2022), 2456-2468.
\bibitem{CommutativeRingTheory} { H. Matsumura}, \textsl{Commutative Ring Theory},  Cambridge Univ., 1987.
\bibitem{VanishingOfCohomologyOverCompleteIntersectionRings} { A. Sadeghi}, Vanishing of cohomology over complete intersection rings, Glasg. Math. J. \textbf{57} (2015), 445--455. 
\bibitem{SemidualizingModules} { S. Sather-Wagstaff}, \textit{Semidualizing modules}, available at \url{https://ssather.people.clemson.edu/DOCS/sdm.pdf}.
\bibitem{HomologicalAspectsOfSemidualizingModules} { R. Takahashi, D.White},  Homological aspects of semidualizing modules, Math. Scand. \textbf{106} (2010), 5--22.
\bibitem{DivisorTheoryinModuleCategories} { W. Vasconcelos},  \textsl{Divisor Theory in Module Categories}, North-Holland Publishing Co., 
1974.
\bibitem{weibel} { C. Weibel}, \textsl{An Introduction to Homological Algebra}, Cambridge Univ. Press, 1994.
\bibitem{livroMCMyoshino} { Y. Yoshino}, \textsl{Maximal Cohen-Macaulay Modules over Cohen-Macaulay Rings}, Cambridge Univ. Press, 1990. 
\bibitem{NumericalApectsofComplexesOfFiniteHomologicalDimensions} { M. Zargar,  M. Gheibi}, Numerical aspects of complexes of finite homological dimensions, arXiv preprint arXiv: 2305.12251v1, 2023.
\end{thebibliography}
     \end{document}